\documentclass[twoside, 12pt]{article}
\usepackage{amssymb,amsmath,amsfonts,soul,amsthm,enumerate}
\usepackage{amsfonts}
\usepackage{graphicx}
\usepackage[cp1251]{inputenc}
\usepackage[T2A]{fontenc}
\usepackage{mathrsfs}
\usepackage[english]{babel}
 \newcommand{\grant}[1]{\medskip \baselineskip 10pt{\footnotesize #1} \medskip}

\newtheorem{theorem}{Theorem}[section]

\newtheorem{remark}{Remark}

\numberwithin{equation}{section}
 \def\@evenhead{\vbox{\hbox to \textwidth{\thepage\hfil\sl\leftmark\strut}\hrule}}
 \def\@oddhead{\vbox{\hbox to \textwidth{\rightmark\hfill\thepage\strut}\hrule}}

\begin{document}
 \sloppy

\centerline{\bf UNIQUENESS OF AN INVERSE SOURCE}
\vskip 0.1cm

\centerline{\bf NON-LOCAL PROBLEM FOR}
\vskip 0.1cm

\centerline{\bf FRACTIONAL ORDER MIXED TYPE EQUATION}     

\vskip 0.3cm

\centerline{\bf M.~S.~Salakhitdinov, E.~T.~Karimov}        

\markboth{\hfill{\footnotesize\rm   M.~S.~Salakhitdinov, E.~T.~Karimov  }\hfill}
{\hfill{\footnotesize\sl  Uniqueness of an inverse source non-local problem for fractional order mixed type equation}\hfill}
\vskip 0.3cm

\vskip 0.7 cm


\vskip 0.2cm


\vskip 0.2cm

\noindent {\bf Abstract.}

In the present work, we investigate a uniqueness of solution of the inverse source problem with non-local conditions for mixed parabolic-hyperbolic type equation with Caputo fractional derivative. Solution of the problem we represent as bi-orthogonal series with respect to space variable and will get fractional order differential equations with respect to time-variable. Using boundary and gluing conditions, we deduce system of algebraic equations regarding unknown constants and imposing condition to the determinant of this system, we prove a uniqueness of considered problem. Moreover, we find some non-trivial solutions of the problem in case, when imposed conditions are not fulfilled.

\section{\large Formulation of a problem}

Consider an equation
\begin{equation}\label{eq.1.1}
f(x)=\left\{
\begin{aligned}
&_CD_{0t}^\alpha u-u_{xx},\,\,\,\,\,\,t> 0,\\
&_CD_{t0}^\beta u-u_{xx},\,\,\,\,\,\,t<0
\end{aligned}
\right.
\end{equation}
in a rectangular domain $\Omega=\left\{(x,t):\, 0<x<1,\, -p<t<q\right\}$. Here $\alpha,\beta, p, q \in \mathbb{R}$ such that $0<\alpha\le1,\, 1<\beta\le2$, $f(x)$ is unknown function,
$$
_{C}D_{0t}^{\alpha }g=\left\{ \begin{aligned}
  & \frac{1}{\Gamma \left( 1-\alpha  \right)}\int\limits_{0}^{t}{\frac{{g}'\left( z \right)}{{{\left( t-z \right)}^{\alpha }}}dz,\,}\,\,0<\alpha <1, \\
 & \frac{dg}{dt},\,\,\,\,\,\,\,\,\,\,\,\,\,\,\,\,\,\,\,\,\,\,\,\,\,\,\,\,\,\,\,\,\,\,\,\,\,\,\,\,\,\,\,\,\alpha =1, 
\end{aligned} \right.
$$
$$_{C}D_{t0}^{\beta }g=\left\{ \begin{aligned}
  & \frac{1}{\Gamma \left( 2-\beta  \right)}\int\limits_{t}^{0}{\frac{{g}''\left( z \right)}{{{\left( z-t \right)}^{\beta -1}}}dz,\,}\,\,1<\beta <2, \\
 & \frac{{{d}^{2}}g}{d{{t}^{2}}},\,\,\,\,\,\,\,\,\,\,\,\,\,\,\,\,\,\,\,\,\,\,\,\,\,\,\,\,\,\,\,\,\,\,\,\,\,\,\,\,\,\,\,\,\beta =2 
\end{aligned} \right.
$$
are Caputo fractional differential operators [1, 92 p., form. (2.4.16)].

\textbf{Problem.} Find a pair of functions $\left(u(x,t), f(x)\right)$ in a domain $\Omega$, satisfying

i) regularity conditions $u(x,t)\in C\left(\overline{\Omega}\right), u_{xx}(x,t)\in C^2\left(\Omega^+\cup \Omega^-\right)$, ${}_CD_{0t}^\alpha u \in C\left(\Omega^+\right)$,\\
 $\,\,\,\,\,\,{}_CD_{t0}^\beta u \in C(\Omega^-), \,f(t)\in C(0,1)$;\\

ii) equation (1.1) in $\Omega^+$, $\Omega^-$;

iii) boundary conditions
\begin{equation}\label{eq.1.2}
u(0,t)=u(1,t),\,\,u_x(0,t)=0,\,\,\,-p\leq t\leq q,
\end{equation}
\begin{equation}\label{eq.1.3}
u(x,-p)=0,\,\,u(x,q)=0,\,\,0\leq x\leq 1,
\end{equation}

iv) and transmitting condition
\begin{equation}\label{eq.1.4}
\lim\limits_{t\rightarrow +0} {}_CD_{0t}^\alpha u(x,t)=\lim\limits_{t\rightarrow -0}\frac{\partial u(x,t)}{\partial (-t)},\,\,0<x<1,
\end{equation}
where $\Omega^+=\Omega\cap\{t>0\},\,$ $\Omega^-=\Omega\cap\{t<0\}$.

Solution of this problem we represent as follows:
\begin{equation}\label{eq.1.5}
u(x,t)=V_0(t)+\sum\limits_{k=1}^\infty V_{1k}(t)\cos 2k\pi x+\sum\limits_{k=1}^\infty V_{2k}(t)\cdot x \sin 2k\pi x,\,\,t\geq0,
\end{equation}
\begin{equation}\label{eq.1.6}
u(x,t)=W_0(t)+\sum\limits_{k=1}^\infty W_{1k}(t)\cos 2k\pi x+\sum\limits_{k=1}^\infty W_{2k}(t)\cdot x \sin 2k\pi x,\,\,t\leq0,
\end{equation}
\begin{equation}\label{eq.1.7}
f(x)=f_0+\sum\limits_{k=1}^\infty f_{1k}\cos 2k\pi x+\sum\limits_{k=1}^\infty f_{2k} \cdot x \sin 2k\pi x,
\end{equation}
where
\begin{equation}\label{eq.1.8}
\begin{aligned}
&V_0(t)=2\int\limits_0^1 u(x,t)(1-x)\,dx,\,t\geq 0,\\
&V_{1k}(t)=4\int\limits_0^1 u(x,t) (1-x)\cos 2k\pi x \,dx,\,t\geq 0,\\
&V_{2k}(t)=4\int\limits_0^1 u(x,t)\sin 2k\pi x \,dx,\, t\geq 0,\\
&W_0(t)=2\int\limits_0^1 u(x,t)(1-x)\,dx,\, t\leq 0,\\
&W_{1k}(t)=4\int\limits_0^1 u(x,t) (1-x)\cos 2k\pi x \,dx,\, t\leq 0,\\
&W_{2k}(t)=4\int\limits_0^1 u(x,t)\sin 2k\pi x\, dx,\, t\leq 0,\\
&f_0=2\int\limits_0^1 f(x)(1-x)\,dx,\\
&f_{1k}=4\int\limits_0^1 f(x) (1-x)\cos 2k\pi x\, dx,\\
&f_{2k}=4\int\limits_0^1 f(x)\sin 2k\pi x\, dx.
\end{aligned}
\end{equation}

Detailed explanation of this representation can be found in [2, p.62], which is based on [3,4].

We would like note some works [5-7], where local and non-local inverse source problems for time-fractional diffusion and diffusion-wave equations were studied. Especially, work by M.Kirane and S.A.Malik [8], where similar non-local conditions were in use.

Based on (1.8), we introduce similar functions with small shift into to the interior of the considered domain. Then applying appropriate Caputo fractional operators, after integrating by parts, we deduce
\begin{equation}\label{eq.1.9}
{}_CD_{0t}^\alpha V_0(t)=f_0,\,\,\,\,\, t\ge0,
\end{equation}
\begin{equation}\label{eq.1.10}
{}_CD_{t0}^\beta W_0(t)=f_0,\,\,\,\,\,\, t<0,
\end{equation}
\begin{equation}\label{eq.1.11}
{}_CD_{0t}^\alpha V_{1k}(t)+(2k\pi)^2V_{1k}(t)=f_{1k}+4k\pi V_{2k}(t),\,\,\,\,\,\, t\ge0,
\end{equation}
\begin{equation}\label{eq.1.12}
{}_CD_{t0}^\beta W_{1k}(t)+(2k\pi)^2W_{1k}(t)=f_{1k}+4k\pi W_{2k}(t),\,\,\,\,\,\, t<0,
\end{equation}
\begin{equation}\label{eq.1.13}
{}_CD_{0t}^\alpha V_{2k}(t)+(2k\pi)^2V_{2k}(t)=f_{2k},\,\,\,\,\,\, t\ge0,
\end{equation}
\begin{equation}\label{eq.1.14}
{}_CD_{t0}^\beta W_{2k}(t)+(2k\pi)^2W_{2k}(t)=f_{2k},\,\,\,\,\,\, t<0.
\end{equation}
General solutions of (1.9) and (1.13) can be written as [1, p.231, form. (4.1.66)]
\begin{equation}\label{eq.1.15}
V_0(t)=V_0(0)+\frac{f_0}{\Gamma(\alpha+1)}t^\alpha,
\end{equation}
\begin{equation}\label{eq.1.16}
V_{2k}(t)=V_{2k}(0)E_{\alpha,1}\left(-(2k\pi)^2t^\alpha\right)+f_{2k}t^\alpha E_{\alpha,\alpha+1}\left(-(2k\pi)^2t^\alpha\right),
\end{equation}
respectively. Here
$$
E_{\alpha,\beta}(t)=\sum\limits_{k=0}^\infty \frac{z^k}{\Gamma(\alpha k +\beta)},\,\,\,\, \alpha,\beta>0
$$
is the Mittag-Leffler function of two parameters [9, p.17].

General solution of equation (1.11) we write as
\begin{equation}\label{eq.1.17}
\begin{aligned}
  & {{V}_{1k}}(t)={{V}_{1k}}(0){{E}_{\alpha ,1}}\left( -{{(2k\pi )}^{2}}{{t}^{\alpha }} \right)+{{f}_{1k}}\cdot {{t}^{\alpha }}\cdot {{E}_{\alpha ,\alpha +1}}\left( -{{(2k\pi )}^{2}}{{t}^{\alpha }} \right)+ \\
 & +4k\pi \cdot {{V}_{2k}}(0)\cdot {{t}^{\alpha }}\cdot {{E}_{1}}\left( \begin{matrix}
   1,1;1,1 & |-{{(2k\pi )}^{2}}{{t}^{\alpha }}  \\
   \alpha +1,\alpha ,\alpha ;1,1;1,1 & |-{{(2k\pi )}^{2}}{{t}^{\alpha }}  \\
\end{matrix} \right)+ \\
 & +4k\pi \cdot {{f}_{2k}}\cdot {{t}^{2\alpha }}\cdot {{E}_{1}}\left( \begin{matrix}
   1,1;1,1 & |-{{(2k\pi )}^{2}}{{t}^{\alpha }}  \\
   2\alpha +1,\alpha ,\alpha ;1,1;1,1 & |-{{(2k\pi )}^{2}}{{t}^{\alpha }}  \\
\end{matrix} \right). 
\end{aligned}
\end{equation}
Based on [1, p.232, form. (4.1.74)], general solution of (1.10) we write as
\begin{equation}\label{eq.1.18}
W_0(t)=W_0(0)-tW_0'(0)+\frac{f_0}{\Gamma(\beta+1)} (-t)^\beta.
\end{equation}
Similarly, we can write general solution of (1.14) as follows
\begin{equation}\label{eq.1.19}
\begin{aligned}
&W_{2k}(t)=W_{2k}(0)E_{\beta,1} \left(-(2k\pi)^2(-t)^\beta\right)-tW_{2k}'(0)E_{\beta,2} \left(-(2k\pi)^2(-t)^\beta\right)+\\
&+f_{2k} (-t)^\beta E_{\beta,\beta+1} \left(-(2k\pi)^2(-t)^\beta\right).
\end{aligned}
\end{equation}
General solution of (1.12) has a form
\begin{equation}\label{eq.1.20}
\begin{aligned}
  & {{W}_{1k}}(t)={{W}_{1k}}(0){{E}_{\beta ,1}}\left( -{{(2k\pi )}^{2}}{{(-t)}^{\beta }} \right)-t{{W}_{1k}}^{\prime }(0){{E}_{\beta ,2}}\left( -{{(2k\pi )}^{2}}{{(-t)}^{\beta }} \right)+ \\
 & +{{f}_{1k}}\cdot {{(-t)}^{\beta }}{{E}_{\beta ,\beta +1}}\left( -{{(2k\pi )}^{2}}{{(-t)}^{\beta }} \right)+ \\
 & +4k\pi \cdot {{W}_{2k}}(0)\cdot {{(-t)}^{\beta }}{{E}_{1}}\left( \begin{matrix}
   1,1;1,1 & |-{{(2k\pi )}^{2}}{{\left( -t \right)}^{\beta }}  \\
   \beta +1,\beta ,\beta ;1,1;1,1 & |-{{(2k\pi )}^{2}}{{\left( -t \right)}^{\beta }}  \\
\end{matrix} \right)+ \\
 & +4k\pi \cdot {{W}_{2k}}^{\prime }(0)\cdot {{(-t)}^{\beta +1}}{{E}_{1}}\left( \begin{matrix}
   1,1;1,1 & |-{{(2k\pi )}^{2}}{{\left( -t \right)}^{\beta }}  \\
   \beta +2,\beta ,\beta ;1,1;1,1 & |-{{(2k\pi )}^{2}}{{\left( -t \right)}^{\beta }}  \\
\end{matrix} \right)+ \\
 & +4k\pi \cdot {{f}_{2k}}\cdot {{(-t)}^{2\beta }}{{E}_{1}}\left( \begin{matrix}
   1,1;1,1 & |-{{(2k\pi )}^{2}}{{\left( -t \right)}^{\beta }}  \\
   2\beta +1,\beta ,\beta ;1,1;1,1 & |-{{(2k\pi )}^{2}}{{\left( -t \right)}^{\beta }}  \\
\end{matrix} \right), 
\end{aligned}
\end{equation}
Here
$$
E_1
\left(
\begin{array}{ll}
\gamma_1, \alpha_1; \gamma_2, \beta_1 &|  x\\
\delta_1, \alpha_2, \beta_2; \delta_2, \alpha_3; \delta_3, \beta_3 &|  y\\
                         \end{array}
                       \right)=\sum\limits_{m,n=0}^\infty \frac{(\gamma_1)_{\alpha_1 m}(\gamma_2)_{\beta_1 n}}{\Gamma(\delta_1+\alpha_2 m+\beta_2 n)}\cdot \frac{x^m}{\Gamma(\delta_2+\alpha_3 m)}\cdot \frac{y^n}{\Gamma(\delta_3+\beta_3 n)}
$$
is the Mittag-Leffler type function in two variables, introduced by Garg et al in [10, form. (11)].

Now, using conditions (1.2)-(1.4) we find unknown constants $f_0,f_{1k}, f_{2k}$, $V_0(0), V_{1k}(0), V_{2k}(0)$, $W_0(0), W_{1k}(0),W_{2k}(0), W_0'(0), W_{1k}'(0), W_{2k}'(0)$.

Found solutions we substitute into the condition (1.3) and deduce
\begin{equation}\label{eq.1.21}
{{W}_{0}}(0)+p{{W}_{0}}^{\prime }(0)+\frac{{{f}_{0}}}{\Gamma \left( \beta +1 \right)}{{p}^{\beta }}=0,
\end{equation}
\begin{equation}\label{eq.1.22}
\begin{aligned}
  & {{W}_{1k}}(0){{E}_{\beta ,1}}\left( -{{(2k\pi )}^{2}}{{p}^{\beta }} \right)+p{{W}_{1k}}^{\prime }(0){{E}_{\beta ,2}}\left( -{{(2k\pi )}^{2}}{{p}^{\beta }} \right)+{{f}_{1k}}{{p}^{\beta }}{{E}_{\beta ,\beta +1}}\left( -{{(2k\pi )}^{2}}{{p}^{\beta }} \right)+ \\
 & +4k\pi \cdot {{W}_{2k}}(0)\cdot {{p}^{\beta }}{{E}_{1}}\left( \begin{matrix}
   1,1;1,1 & |-{{(2k\pi )}^{2}}{{p}^{\beta }}  \\
   \beta +1,\beta ,\beta ;1,1;1,1 & |-{{(2k\pi )}^{2}}{{p}^{\beta }}  \\
\end{matrix} \right)+ \\
 & +4k\pi \cdot {{W}_{2k}}^{\prime }(0)\cdot {{p}^{\beta +1}}{{E}_{1}}\left( \begin{matrix}
   1,1;1,1 & |-{{(2k\pi )}^{2}}{{p}^{\beta }}  \\
   \beta +2,\beta ,\beta ;1,1;1,1 & |-{{(2k\pi )}^{2}}{{p}^{\beta }}  \\
\end{matrix} \right)+ \\
 & +4k\pi \cdot {{f}_{2k}}\cdot {{p}^{2\beta }}{{E}_{1}}\left( \begin{matrix}
   1,1;1,1 & |-{{(2k\pi )}^{2}}{{p}^{\beta }}  \\
   2\beta +1,\beta ,\beta ;1,1;1,1 & |-{{(2k\pi )}^{2}}{{p}^{\beta }}  \\
\end{matrix} \right)=0, 
\end{aligned}
\end{equation}
\begin{equation}\label{eq.1.23}
\begin{aligned}
& {{W}_{2k}}(0){{E}_{\beta ,1}}\left( -{{(2k\pi )}^{2}}{{p}^{\beta }} \right)+p{{W}_{2k}}'(0){{E}_{\beta ,2}}\left( -{{(2k\pi )}^{2}}{{p}^{\beta }} \right)+\\
& +{{f}_{2k}}{{p}^{\beta }}{{E}_{\beta ,\beta +1}}\left( -{{(2k\pi )}^{2}}{{p}^{\beta }} \right)=0,
\end{aligned}
\end{equation}
\begin{equation}\label{eq.1.24}
{{V}_{0}}(0)+\frac{{{f}_{0}}}{\Gamma \left( \alpha +1 \right)}{{q}^{\alpha }}=0,
\end{equation}
\begin{equation}\label{eq.1.25}
\begin{aligned}
  & {{V}_{1k}}(0){{E}_{\alpha ,1}}\left( -{{(2k\pi )}^{2}}{{q}^{\alpha }} \right)+{{f}_{1k}}{{q}^{\alpha }}{{E}_{\alpha ,\alpha +1}}\left( -{{(2k\pi )}^{2}}{{q}^{\alpha }} \right)+ \\
 & +4k\pi \cdot {{V}_{2k}}(0)\cdot {{q}^{\alpha }}{{E}_{1}}\left( \begin{matrix}
   1,1;1,1 & |-{{(2k\pi )}^{2}}{{q}^{\alpha }}  \\
   \alpha +1,\alpha ,\alpha ;1,1;1,1 & |-{{(2k\pi )}^{2}}{{q}^{\alpha }}  \\
\end{matrix} \right)+ \\
 & +4k\pi \cdot {{f}_{2k}}\cdot {{q}^{2\alpha }}{{E}_{1}}\left( \begin{matrix}
   1,1;1,1 & |-{{(2k\pi )}^{2}}{{q}^{\alpha }}  \\
   2\alpha +1,\alpha ,\alpha ;1,1;1,1 & |-{{(2k\pi )}^{2}}{{q}^{\alpha }}  \\
\end{matrix} \right)=0, 
\end{aligned}
\end{equation}
\begin{equation}\label{eq.1.26}
{{V}_{2k}}(0){{E}_{\alpha ,1}}\left( -{{(2k\pi )}^{2}}{{q}^{\alpha }} \right)+{{f}_{2k}}{{q}^{\alpha }}{{E}_{\alpha ,\alpha +1}}\left( -{{(2k\pi )}^{2}}{{q}^{\alpha }} \right)=0.
\end{equation}
Considering $u(x,+0)=u(x,-0)$, which follows from $u(x,t)\in C\left(\overline{\Omega}\right)$, we get
\begin{equation}\label{eq.1.27}
V_0(0)=W_0(0),\,\,V_{1k}(0)=W_{1k}(0),\,\,V_{2k}(0)=W_{2k}(0).
\end{equation}
Based on (1.9), (1.11), (1.13), we deduce
\begin{equation}\label{eq.1.28}
\begin{aligned}
&\lim\limits_{t\rightarrow +0} {}_CD_{0t}^\alpha V_0(t)=f_0,\\
&\lim\limits_{t\rightarrow +0} {}_CD_{0t}^\alpha V_{1k}(t)=f_{1k}+4k\pi V_{2k}(0)-{{\left( 2k\pi  \right)}^{2}}{{V}_{1k}}\left( 0 \right),\\
&\lim\limits_{t\rightarrow +0} {}_CD_{0t}^\alpha V_{2k}(t)=f_{2k}-(2k\pi)^2V_{2k}(0).
\end{aligned}
\end{equation}
Now we calculate $\frac{d}{d(-t)}W_0(t), \frac{d}{d(-t)}W_{1k}(t)$ and $\frac{d}{d(-t)}W_{2k}(t)$. One can easily deduce that
\begin{equation}\label{eq.1.29}
\frac{d}{d(-t)}W_0(t)=W_0'(0)+\frac{f_0}{\Gamma(\beta)}(-t)^{\beta-1}.
\end{equation}
From (1.19) we find
$$
\begin{aligned}
&\frac{d}{d(-t)}W_{2k}(t)=W_{2k}(0)\frac{d}{d(-t)}E_{\beta,1}\left(-(2k\pi)^2 (-t)^\beta\right)+\\
&+W_{2k}'(0)\frac{d}{d(-t)}\left(-tE_{\beta,2}\left(-(2k\pi)^2 (-t)^\beta\right)\right)+\\
&+f_{2k}\frac{d}{d(-t)}\left((-t)^\beta E_{\beta,\beta+1}\left(-(2k\pi)^2 (-t)^\beta\right)\right).
\end{aligned}
$$
We use formula of differentiation (see [9], p.21, form. (1.82))
$$
{}_{RL}D_{0t}^\gamma\left(t^{\alpha k+\beta-1} E_{\alpha,\beta}^{(k)}\left(\lambda t^\alpha\right)\right)=t^{\alpha k+\beta-\gamma-1} E_{\alpha,\beta-\gamma}\left((\lambda t^\alpha\right),
$$
where ${}_{RL}D_{0t}^\gamma (\cdot)$ is the Riemann-Liouville fractional derivative of order $\gamma$ [1], $E_{\alpha,\beta}^{(k)}(t)=\frac{d^k}{dt^k}E_{\alpha,\beta}(t)$. After some evaluations we deduce
\begin{equation}\label{eq.1.30}
\begin{aligned}
&\frac{d}{d(-t)}W_{2k}(t)=-W_{2k}(0)(2k\pi)^2 (-t)^{\beta-1} E_{\beta,\beta}\left(-(2k\pi)^2 (-t)^\beta\right)+\\
&+W_{2k}'(0)E_{\beta,1}\left(-(2k\pi)^2 (-t)^\beta\right)+\\
&+f_{2k}(-t)^{\beta-1} E_{\beta,\beta}\left(-(2k\pi)^2 (-t)^\beta\right).
\end{aligned}
\end{equation}
By similar way, using the following formula (see [10], form. (33))
$$
\begin{aligned}
&{}_{RL}D_{ax}^\gamma\left\{(x-a)^{\delta_1-1}E_1
\left(
\begin{array}{ll}
\gamma_1, \alpha_1;\gamma_2,\beta_1 &|  w_1(x-a)^{\alpha_2}\\
\delta_1,\alpha_2,\beta_2;\delta_2,\alpha_3,\delta_3,\beta_3 &|  w_2(x-a)^{\beta_2}\\
                         \end{array}\right)
                         \right\}=\\
&(x-a)^{\delta_1-\gamma-1}E_1
\left(
\begin{array}{ll}
\gamma_1, \alpha_1; \gamma_2, \beta_1 &|  w_1(x-a)^{\alpha_2}\\
\delta_1-\gamma, \alpha_2, \beta_2; \delta_2, \alpha_3; \delta_3, \beta_3 &|  w_2(x-a)^{\beta_2}\\
                         \end{array} \right),
\end{aligned}
$$
we deduce
\begin{equation}\label{eq.1.31}
\begin{aligned}
&\frac{d}{d(-t)}W_{1k}(t)=-(2k\pi)^2 W_{1k}(0)(-t)^{\beta-1}E_{\beta,\beta}\left(-(2k\pi)^2 (-t)^\beta\right)+\\
&+W_{1k}'(0) E_{\beta,1}\left(-(2k\pi)^2 (-t)^\beta\right)
+f_{1k} (-t)^{\beta-1} E_{\beta,\beta+1}\left(-(2k\pi)^2 (-t)^\beta\right)+\\
&+4k\pi W_{2k}(0)(-t)^{\beta-1} {{E}_{1}}\left( \begin{matrix}
   1,1;1,1 & |-{{(2k\pi )}^{2}}{{(-t)}^{\beta }}  \\
   \beta,\beta ,\beta ;1,1;1,1 & |-{{(2k\pi )}^{2}}{{(-t)}^{\beta}}  \\
\end{matrix} \right)+\\
&+4k\pi W_{2k}'(0) (-t)^\beta {{E}_{1}}\left( \begin{matrix}
   1,1;1,1 & |-{{(2k\pi )}^{2}}{{(-t)}^{\beta }}  \\
   \beta+1,\beta ,\beta ;1,1;1,1 & |-{{(2k\pi )}^{2}}{{(-t)}^{\beta}}  \\
\end{matrix} \right)+\\
&+4k\pi f_{2k} (-t)^{2\beta-1}{{E}_{1}}\left( \begin{matrix}
   1,1;1,1 & |-{{(2k\pi )}^{2}}{{(-t)}^{\beta }}  \\
   2\beta,\beta ,\beta ;1,1;1,1 & |-{{(2k\pi )}^{2}}{{(-t)}^{\beta}}  \\
\end{matrix} \right).
\end{aligned}
\end{equation}
From (1.29)-(1.31) we get
\begin{equation}\label{eq.1.32}
\lim\limits_{t\rightarrow -0}\frac{d}{d(-t)}W_0(t)=W_0'(0),\,\,\lim\limits_{t\rightarrow -0} \frac{d}{d(-t)}W_{1k}(t)=W_{1k}'(0),\,\lim\limits_{t\rightarrow -0} \frac{d}{d(-t)}W_{2k}(t)=W_{2k}'(0).
\end{equation}
Considering (1.28) and (1.32) we have
\begin{equation}\label{eq.1.33}
\begin{aligned}
&f_0=W_0'(0),\\
&f_{1k}+4k\pi V_{2k}(0)-(2k\pi)^2 V_{1k}(0)=W_{1k}'(0),\\
&f_{2k}-(2k\pi)^2 V_{2k}(0)=W_{2k}'(0).
\end{aligned}
\end{equation}
From (1.21), (1.24) and first relations of (1.27), (1.33) we get the following system of equations
\begin{equation}\label{eq.1.34}
\left\{
\begin{aligned}
&W_0(0)+\left(\frac{p^\beta}{\Gamma(\beta+1)}+p\right)W_0'(0)=0,\\
&W_0(0)+\frac{q^\alpha}{\Gamma(\alpha+1)}W_0'(0)=0,\\
&f_0=W_0'(0).
\end{aligned}
\right.
\end{equation}
If
\begin{equation}\label{eq.1.35}
\Delta_0=p+\frac{p^\beta}{\Gamma(\beta+1)}-\frac{q^\alpha}{\Gamma(\alpha+1)}\neq 0,
\end{equation}
then we find that
\begin{equation}\label{eq.1.36}
f_0=W_0'(0)=V_0(0)=W_0(0)=0.
\end{equation}

Now from (1.22), (1.23), (1.25), (1.26) and last two relations of (1.27), (1.33), we obtain another system of equations
\begin{equation}\label{eq.1.37}
\left\{
\begin{aligned}
&f_{1k}=W_{1k}'(0)+(2k\pi)^2W_{1k}(0)-4k\pi W_{2k}(0),\\
&f_{2k}=W_{2k}'(0)+(2k\pi)^2 W_{2k}(0),\\
&W_{2k}(0)+W_{2k}'(0)\left[p^\beta E_{\beta,\beta+1}\left(-(2k\pi)^2 p^\beta\right)+p E_{\beta,2}\left(-(2k\pi)^2 p^\beta\right)\right]=0,\\
&W_{2k}(0)+W_{2k}'(0) q^\alpha E_{\alpha,\alpha+1}\left(-(2k\pi)^2 q^\alpha\right)=0,\\
&W_{1k}(0)+W_{1k}'(0)\left[p^\beta E_{\beta,\beta+1}\left(-(2k\pi)^2 p^\beta\right)+p E_{\beta,2}\left(-(2k\pi)^2 p^\beta\right)\right]=0,\\
&W_{1k}(0)+W_{1k}'(0)q^\alpha E_{\alpha,\alpha+1}\left(-(2k\pi)^2 q^\alpha\right)=0.
\end{aligned}
\right.
\end{equation}
Further, assuming
\begin{equation}\label{eq.1.38}
\Delta_k=p^\beta E_{\beta,\beta+1}\left(-(2k\pi)^2 p^\beta\right)+p E_{\beta,2}\left(-(2k\pi)^2 p^\beta\right)-q^\alpha E_{\alpha,\alpha+1}\left(-(2k\pi)^2 q^\alpha\right)\neq 0,
\end{equation}
we get
\begin{equation}\label{eq.1.39}
V_{1k}(0)=V_{2k}(0)=W_{1k}(0)=W_{2k}(0)=W_{1k}'(0)=W_{2k}'(0)=f_{1k}=f_{2k}=0.
\end{equation}

If conditions (1.35), (1.38) hold, then based on (1.5)-(1.7), due to (1.36), (1.39) we have
$$
\begin{aligned}
&\int\limits_0^1 u(x,t)(1-x)\,dx=0,\,\, \int\limits_0^1 u(x,t)\sin 2k\pi x\, dx=0,\\
&\int\limits_0^1 u(x,t)(1-x)\cos 2k\pi x dx=0,\,\,\int\limits_0^1 f(x) (1-x)\,dx=0,\\
&\int\limits_0^1 f(x)\sin 2k\pi x \,dx=0,\,\, \int\limits_0^1 f(x) (1-x) \cos 2k\pi x\, dx=0,\, k=1,2,...
\end{aligned}
$$

According to the completeness of the system $\left\{ 1, \cos 2k\pi x, x\sin 2k\pi x\right\}$ in $L_2[0,1]$, we can state that $u(x,t)=0$ a.e. in $[0,1]$ for $t\in [-p,q]$ and $f(x)=0$ a.e. in $[0,1]$.

Now we formulated above obtained result as the following
\begin{theorem}
If conditions (1.35), (1.38) hold, then problem has only trivial solution.
\end{theorem}
\section{Nontrivial solutions of the problem}
We consider case, when conditions (1.35), (1.38) are not fulfilled. Let $\Delta_0=0$ for some $p,\, q$. Then problem has nontrivial solution of the form
\begin{equation}\label{eq.2.1}
u(x,t)=u_0(t),\,\,f(x)=f_0,
\end{equation}
here
$$
u_0(t)=\left\{
\begin{aligned}
&\frac{t^\alpha-q^\alpha}{\Gamma(\alpha+1)}f_0,\,\,t\geq 0,\\
&\left[\frac{(-t)^\beta-p^\beta}{\Gamma(\beta+1)}-t-p\right]f_0,\,\,t\leq 0,
\end{aligned}
\right.
$$
$f_0\neq 0$ is arbitrary constant.

If $\Delta_k=0$ for $k=m\in \mathbb{N}$, i.e. $\Delta_m=0$, then considered problem has nontrivial solutions of the form
\begin{equation}\label{eq.2.2}
u_m(x,t)=\left\{
\begin{aligned}
&V_{1m}(t)\cos 2m\pi x+V_{2m}(t) x \sin 2m \pi x,\,\,t\geq 0,\\
&W_{1m}(t)\cos 2m\pi x+W_{2m}(t) x \sin 2m \pi x,\,\,t\leq 0,\\
\end{aligned}
\right.
\end{equation}
\begin{equation}\label{eq.2.3}
\begin{aligned}
&f_m(x)=\left\{E_{\alpha,1}\left(-(2m\pi)^2 q^\alpha\right)W_{1m}'(0)
+\right.\\
&\left.+4k\pi q^\alpha E_{\alpha,\alpha+1}\left(-(2m\pi)^2 q^\alpha\right)W_{2m}'(0)\right\}\cos2m\pi x+\\
&+E_{\alpha,1}\left(-(2m\pi)^2 q^\alpha\right)W_{2m}'(0) x \sin 2m\pi x,
\end{aligned}
\end{equation}
where
$$
\begin{aligned}
&V_{1m}(t)=W_{1m}'(0)\left\{t^\alpha E_{\alpha,\alpha+1}\left(-(2m\pi)^2 t^\alpha\right)E_{\alpha,1}\left(-(2m\pi)^2 q^\alpha\right)-q^\alpha E_{\alpha,\alpha+1}\left(-(2m\pi)^2 q^\alpha\right)\times \right.\\
&-\left.\times E_{\alpha,1}\left(-(2m\pi)^2 t^\alpha\right)\right\}+
4k\pi t^\alpha W_{2m}'(0)\left\{{}_{}^{}q^\alpha E_{\alpha,\alpha+1}\left(-(2m\pi)^2 q^\alpha\right)\times \right.\\
&\times \left[E_{\alpha,\alpha+1}\left(-(2m\pi)^2 t^\alpha\right)-
{{E}_{1}}\left( \begin{matrix}
   1,1;1,1 & |-{{(2k\pi )}^{2}}{{t}^{\alpha }}  \\
   \alpha +1,\alpha ,\alpha ;1,1;1,1 & |-{{(2k\pi )}^{2}}{{t}^{\alpha }}  \\
\end{matrix} \right)\right]+\\
&\left.+t^\alpha E_{\alpha,1}\left(-(2m\pi)^2 q^\alpha\right)
{{E}_{1}}\left( \begin{matrix}
   1,1;1,1 & |-{{(2k\pi )}^{2}}{{t}^{\alpha }}  \\
   2\alpha +1,\alpha ,\alpha ;1,1;1,1 & |-{{(2k\pi )}^{2}}{{t}^{\alpha }}  \\
\end{matrix} \right)\right\},
\end{aligned}
$$
$$
\begin{aligned}
&V_{2m}(t)=\left\{t^\alpha E_{\alpha,\alpha+1}\left(-(2m\pi)^2 t^\alpha\right)E_{\alpha,1}\left(-(2m\pi)^2 q^\alpha\right)-\right.\\
&\left.-q^\alpha E_{\alpha,\alpha+1}\left(-(2m\pi)^2 q^\alpha\right)E_{\alpha,1}\left(-(2m\pi)^2 t^\alpha\right)\right\}W_{2m}'(0),
\end{aligned}
$$
$$
\begin{aligned}
&W_{1m}(t)=\left\{(-t)^\beta E_{\alpha,1}\left(-(2m\pi)^2 q^\alpha\right)E_{\beta,\beta+1}\left(-(2m\pi)^2 (-t)^\beta\right)-\right.\\
&\left.-q^\alpha E_{\alpha,\alpha+1}\left(-(2m\pi)^2 q^\alpha\right)E_{\beta,1}\left(-(2m\pi)^2 (-t)^\beta\right)\right\}W_{1m}'(0)+\\
&+4k\pi (-t)^\beta W_{2m}'(0)\left\{{}_{}^{}q^\alpha E_{\alpha,\alpha+1}\left(-(2m\pi)^2 q^\alpha\right)\times \right.\\
&\times \left[E_{\beta,\beta+1}\left(-(2m\pi)^2 (-t)^\beta\right)-
{{E}_{1}}\left( \begin{matrix}
   1,1;1,1 & |-{{(2k\pi )}^{2}}{(-t)^{\beta }}  \\
   \beta +1,\beta ,\beta ;1,1;1,1 & |-{{(2k\pi )}^{2}}{(-t)^{\beta }}  \\
\end{matrix} \right)\right]-\\
&-t{{E}_{1}}\left( \begin{matrix}
   1,1;1,1 & |-{{(2k\pi )}^{2}}{(-t)^{\beta }}  \\
   \beta +2,\beta ,\beta ;1,1;1,1 & |-{{(2k\pi )}^{2}}{(-t)^{\beta }}  \\
\end{matrix} \right)+\\
&\left.+(-t)^\beta E_{\alpha,1}\left(-(2m\pi)^2 q^\alpha\right){{E}_{1}}\left( \begin{matrix}
   1,1;1,1 & |-{{(2k\pi )}^{2}}{(-t)^{\beta }}  \\
   2\beta +1,\beta ,\beta ;1,1;1,1 & |-{{(2k\pi )}^{2}}{(-t)^{\beta }} \\
\end{matrix} \right)\right\}
\end{aligned}
$$
$$
\begin{aligned}
&W_{2m}(t)=\left\{(-t)^\beta E_{\alpha,1}\left(-(2m\pi)^2 q^\alpha\right)E_{\beta,\beta+1}\left(-(2m\pi)^2 (-t)^\beta\right)-\right.\\ &-q^\alpha E_{\alpha,\alpha+1}\left(-(2m\pi)^2 q^\alpha\right)E_{\beta,1}\left(-(2m\pi)^2 (-t)^\beta\right)-
\left.-tE_{\beta,2}\left(-(2m\pi)^2 (-t)^\beta\right)\right\}W_{2m}'(0).
\end{aligned}
$$
Here $W_{1m}'(0), W_{2m}(0)$ are arbitrary non-zero constants.

\begin{remark}
If conditions (1.35), (1.38) are not fulfilled, there exist nontrivial solutions of the considered problem and they have form (2.1) or (2.2)-(2.3).
\end{remark}

\section*{\large Acknowledgments} The authors thank the unknown referee for his useful suggestions.

 \grant{This work was partially supported by the Grant no. ~3293/GF4 of the Ministry
of education and science of the Republic of Kazakhstan}

\def\bibname{\vspace*{-30mm}{

\vskip 1 cm \footnotesize
\begin{flushleft}
Makhmud Salakhitdinov, Erkinjon Karimov \\ 
Department of Differential Equations\\ 
Institute of Mathematics, National University of Uzbekistan \\ 
29 Durmon yuli St,\\ 
100125 Tashkent, Uzbekistan \\ 
E-mails: salakhitdinovms@yahoo.com, erkinjon@gmail.com
\end{flushleft}

\vskip0.5cm
\begin{flushright}
Received: ??.??.20??
\end{flushright}

 \end{document}